\documentclass[]{article}  
\usepackage{CJK}
\usepackage{amssymb}
\usepackage{bbm}
\usepackage{amsfonts}
\usepackage{mathrsfs}
\usepackage{graphicx}
\usepackage{amsmath}
\usepackage{amsthm}
\usepackage{xypic}
\usepackage{graphicx}
\usepackage{times}
\usepackage{geometry}
\usepackage{color}
\usepackage{indentfirst}
\usepackage{cases}
\usepackage{hyperref}
\usepackage{graphicx}
\usepackage{float}
\usepackage{amsthm}
\usepackage{amssymb}
\usepackage{amsmath}
\usepackage{mathrsfs}
\usepackage{indentfirst}
\usepackage{geometry}
\usepackage{authblk}
\usepackage{hyperref}
\usepackage{enumerate}

\newtheorem{Theorem}{Theorem}[section]
\theoremstyle{definition}
\newtheorem{Corollary}[Theorem]{Corollary}
\newtheorem{Problem}[Theorem]{Problrm}
\newtheorem{Lemma}[Theorem]{Lemma}
\newtheorem{Proposition}[Theorem]{Proposition}

\theoremstyle{definition}

\title{The m-partite digraphical  representations of valency 3 of finite groups generated by two elements}
\author{Songnian Xu \thanks{Corresponding author. E-mail address: xsn131819@163.com},\ \ \ Dein Wong \thanks{Corresponding author. E-mail address:wongdein@163.com.  Supported by the National Natural Science Foundation of China(No.12371025)} ,\ \ Chi Zhang \thanks{ Supported by NSFC of China(No.12001526)}, \ \ Wenhao Zhen
\\ {\small  \it School of Mathematics, China University of Mining and Technology, Xuzhou,  China.}

}
\date{}

\begin{document}
\baselineskip 17pt

\title{The $m$-partite digraphical  representations of valency 3 of finite groups generated by two elements}

\author{Songnian Xu \thanks{Corresponding author. E-mail address: xsn131819@163.com},\ \ \ Dein Wong \thanks{Corresponding author. E-mail address:wongdein@163.com.  Supported by the National Natural Science Foundation of China(No.12371025)} ,\ \ Chi Zhang \thanks{ Supported by NSFC of China(No.12001526)}, \ \ Wenhao Zhen
\\ {\small  \it School of Mathematics, China University of Mining and Technology, Xuzhou,  China.}

}

\date{}
\maketitle

\begin{abstract}
\ \ Let $G$ be a finite group and $m$ be an integer.
We employ the notation $g_i$ to represent elements $(g,i)$ in the Cartesian product $G \times \mathbb{Z}_m$, where $\mathbb{Z}_m$ denotes integers modulo $m$.
For given sets $T_{i,j} \subseteq G$ ($i,j \in \mathbb{Z}_m$), we construct the $m$-$Cayley$ $digraph$ $\Gamma = \mathrm{Cay}(G, T_{i,j}: i,j \in \mathbb{Z}_m)$ with vertex set $\bigcup_{i\in\mathbb{Z}_m}G_i$ (where $G_i = \{g_i | g \in G\}$) and arc set $\bigcup_{i,j}\{(g_i, (tg)_j) | t \in T_{i,j}, g \in G\}$.
When $T_{i,i} = \emptyset$ for all $i \in \mathbb{Z}_m$, we call $\Gamma$ an \emph{$m$-partite Cayley digraph}.
For $m$-partite Cayley digraphs, we observe that a $1$-partite Cayley digraph is necessarily an empty graph. Therefore, throughout this paper, we restrict our consideration to the case where $m \geq 2$.
The digraph $\Sigma$ is regular if there exists a non-negative integer $k$ such that every vertex has out-valency and in-valency equal to $k$. All digraphs considered in this paper are regular.
We say a group $G$ admits an \emph{$m$-partite digraphical representation} ($m$-PDR for short) if there exists a regular $m$-partite Cayley digraph $\Gamma$ with $\mathrm{Aut}(\Gamma) \cong G$.
 Based on Du et al.'s complete classification of unrestricted $m$-PDRs \cite{du4} (2022), we focus on the unresolved valency-specific cases.
 In this paper, we investigate $m$-PDRs of valency 3 for groups generated by at most two elements, and establish a complete classification of nontrivial finite simple groups admitting $m$-PDRs of valency 3 with $m\geq2$.
\end{abstract}


\section{Introduction}

A $digraph$ $\Gamma$ is an ordered pair $(V(\Gamma), A(\Gamma))$, where $V(\Gamma)$ is a non-empty set and $A(\Gamma)$ is a subset of $V(\Gamma) \times V(\Gamma)$. We call $V(\Gamma)$ the \emph{vertex set} and $A(\Gamma)$ the \emph{arc set} of $\Gamma$, and their elements are called \emph{vertices} and \emph{arcs} respectively. For an arc $(u, v) \in A(\Gamma)$, we say $v$ is an \emph{out-neighbor} of $u$ and $u$ is an \emph{in-neighbor} of $v$.
The $out$-$valency$ (resp. $in$-$valency$) of a vertex $v \in V(\Gamma)$ is the number of its out-neighbors (resp. in-neighbors). The digraph $\Gamma$ is \emph{regular} if there exists a non-negative integer $k$ such that every vertex has out-valency and in-valency equal to $k$. All digraphs considered in this paper are regular.
For a subset $X \subseteq V(\Gamma)$, the $subdigraph$ induced by $\Gamma$ on $X$ is $\Gamma[X] := (X, A(\Gamma) \cap (X \times X))$, which we abbreviate as $[X]$. A digraph $\Gamma$ is a $graph$ if $A(\Gamma)$ is symmetric, i.e., $(u,v) \in A(\Gamma)$ implies $(v,u) \in A(\Gamma)$. A digraph is $oriented$ if for any distinct vertices $u,v$, at most one of $(u,v)$ and $(v,u)$ belongs to $A(\Gamma)$.

An $automorphism$ of a digraph $\Gamma = (V(\Gamma), A(\Gamma))$ is a permutation $\sigma$ of $V(\Gamma)$ that preserves the arc relation: $(x^\sigma, y^\sigma) \in A(\Gamma)$ if and only if $(x,y) \in A(\Gamma)$ for all $x,y \in V(\Gamma)$.
Let $G$ be a finite group and $S$ be a subset of $G$. The $Cayley$ $digraph$ $\mathrm{Cay}(G,S)$ has vertex set $G$ and arc set $\{(g,sg) \mid g \in G, s \in S\}$. This is a $Cayley$ $graph$ if and only if $S = S^{-1}$.

For a permutation group $G$ acting on a set $X$, the $stabilizer$ of $x \in X$ is $G_x = \{g \in G \mid x^g = x\}$. The group $G$ is $semiregular$ if $G_x = \{1\}$ for all $x \in X$, and $regular$ if it is both semiregular and transitive.
The concept of Cayley digraphs can be nicely generalized to $m$-$Cayley$ $digraphs$ where regular actions are replaced with semiregular actions. An $m$-Cayley (di)graph $\Gamma$ over a finite group $G$ is defined as a (di)graph which has a semiregular group of automorphisms isomorphic to $G$ with $m$ orbits on its vertex set.

We say a finite group $G$ admits:
\begin{itemize}
\item a graphical regular representation (GRR) if there exists a Cayley graph $\Gamma = \mathrm{Cay}(G,S)$ with $\mathrm{Aut}(\Gamma) \cong G$
\item a digraphical regular representation (DRR) if there exists a Cayley digraph $\Gamma = \mathrm{Cay}(G,S)$ with $\mathrm{Aut}(\Gamma) \cong G$
\item an oriented regular representation (ORR) if there exists an oriented Cayley digraph $\Gamma = \mathrm{Cay}(G,S)$ with $\mathrm{Aut}(\Gamma) \cong G$
\end{itemize}

The GRR and DRR problem asks which groups admit such representations. Babai \cite{bab1} proved that except for $Q_8$, $\mathbb{Z}_2^2$, $\mathbb{Z}_2^3$, $\mathbb{Z}_2^4$, and $\mathbb{Z}_3^2$, every group admits a DRR. Spiga \cite{mor,mor1,spi1} classified groups admitting ORRs. While GRR implies DRR, the converse fails. The GRR turned out to be much more difficult to handle and, after a long series of partial results by various authors \cite{he,im,im1,im2,no,no1,wa}, the classification was completed by Godsil in \cite{god1}.

We say a finite group $G$ admits:
\begin{itemize}
\item a \emph{graphical m-semiregular representation} (GmSR) if there exists a regular $m$-Cayley graph $\Gamma$ over $G$ with $\mathrm{Aut}(\Gamma) \cong G$
\item a \emph{digraphical m-semiregular representation} (DmSR) if there exists a regular $m$-Cayley digraph $\Gamma$ over $G$ with $\mathrm{Aut}(\Gamma) \cong G$
\item an \emph{oriented m-semiregular representation} (OmSR) if there exists a regular oriented $m$-Cayley digraph $\Gamma$ over $G$ with $\mathrm{Aut}(\Gamma) \cong G$
\end{itemize}
Note that G1SR, D1SR and O1SR coincide with GRR, DRR and ORR respectively. Groups admitting GmSR, DmSR or OmSR for all $m$ were classified in \cite{du1,du1'}.

A directed graph is called an \emph{$m$-partite digraph} if its vertex set can be partitioned into $m$ subsets (or parts) such
that no arcs exist within the same part (equivalently, the induced subgraph on each part is an empty digraph).
We now present the formal definition of  \emph{$m$-partite digraphical representation}, which serves as the central concept in this work.

We say a finite group $G$ admits:
\begin{itemize}
\item an \emph{$m$-partite digraphical representation} ($m$-PDR) if there exists a regular $m$-partite Cayley digraph $\Gamma$ with $\mathrm{Aut}(\Gamma) \cong G$
\item an \emph{$m$-partite oriented semiregular representation} ($m$-POSR) if there exists a regular $m$-partite oriented Cayley digraph $\Gamma$ with $\mathrm{Aut}(\Gamma) \cong G$
\end{itemize}
Clearly $m$-POSR implies $m$-PDR, but not conversely.

\textbf{$m$-PDR with prescribed valency $k$}: In contrast to unrestricted GmSR, DmSR, OmSR and $m$-PDR, the classification problem for groups admitting GmSR, DmSR, OmSR or $m$-PDR with prescribed valency remains largely unresolved.
For digraphs, when $k$ is even, the minimal possible value is 2; when $k$ is odd, the minimal value is 1. However, when $k=1$, a connected digraph of valency one is just a directed cycle. Therefore, for even $k$, the smallest interesting case is valency two, while for odd $k$, it is valency three.
Regarding the valency 2 case, Verret and Xia \cite{x2} classified finite simple groups admitting an ORR (i.e., O1SR) of valency 2, proving that every simple group of order at least 5 has an ORR of valency two.
Du et al. \cite{du1} established the complete classification of $m$-PDRs for all finite groups, while their subsequent work \cite{du2} characterized $m$-POSRs of valency 2 specifically for finite groups generated by at most two elements.
These results naturally lead us to investigate the classification problem of $m$-PDRs of valency 3 for finite simple groups.
According to the classification theorem of finite simple groups, we know that all finite simple groups can be generated by at most two elements.
Therefore, in this article we consider $m$-PDRs of valency 3 for groups generated by at most two elements, and we present the following results.

\begin{Theorem}
Let $G = \langle x \rangle \neq 1$ be a cyclic group. Then for every integer $m\geq2$, $G$ admits an $m$-PDR of valency 3, except for when:
\begin{enumerate}
    \item[(i)] $m = 2$ and $o(x) < 5$,
    \item[(ii)] $m = 3$ and $o(x)=2$ .
\end{enumerate}

\end{Theorem}

\begin{Theorem}
Let $G=\langle x,y\rangle$ be a finite simple group. Then $G$ admits a 2-PDR of valency 3.
\end{Theorem}

\begin{Theorem}
Let $G = \langle x, y \rangle$ be a finite  group. Then for every integer $m \geq 3$, $G$ admits an $m$-PDR of valency 3.

\end{Theorem}

The following corollary is an immediate consequence of Theorems~1.1, 1.2 and 1.3.

\begin{Corollary}

Let $G$ be a non-trivial finite simple group. Then $G$ admits an $m$-PDR of valency 3, except for when:
\begin{enumerate}
    \item[(i)] $m = 2$ and $G \cong \mathbb{Z}_2$ or $\mathbb{Z}_3$,
    \item[(ii)] $m = 3$  and $G \cong \mathbb{Z}_2$.
\end{enumerate}

\end{Corollary}

For the case when $G = \mathbb{Z}_1$ (the trivial group), we have not resolved the classification problem of $m$-PDRs of valency 3 for $G$. In fact, this problem is equivalent to the following fundamental question:
For every integer $m$, does there exist a regular digraph $\Gamma$ of order $m$ and valency 3 with trivial automorphism group?

Furthermore, for the case when $m=2$ and $G=\langle x,y\rangle$ is not a finite simple group, the classification problem of $m$-PDRs of valency 3 for $G$ remains unresolved. In our proof attempts, we observed that the difficulty lies in classifying two-generated groups admitting DRRs of valency two.
To conclude this section, we propose the following problems:

\begin{Problem}

\begin{enumerate}

\item For every integer $m$, does there exist a regular digraph $\Gamma$ of order $m$ and valency 3 with trivial automorphism group?
\item Classify finite groups generated by two elements that admit a DRR of valency two.
\item Let $G = \langle x, y \rangle \neq 1$ be a group that is not finite simple. Classify all such groups \( G \) that admit  \( m \)-PDRs of valency 3.
\end{enumerate}
\end{Problem}

\section{Preliminaries and notations}

Let $m$ be a positive integer, and let $G$ be a group. To simplify notation throughout this paper, we denote the element $(g,i)$ of the Cartesian product $G \times \{0,\ldots,m-1\}$ by $g_i$. We frequently identify $\{0,\ldots,m-1\}$ with $\mathbb{Z}_m$, the integers modulo $m$.
Recall that an $m$-Cayley digraph of a finite group $G$ is a digraph admitting a semiregular group of automorphisms isomorphic to $G$ with exactly $m$ orbits on its vertex set. Now, we consider a more concrete definition for $m$-Cayley digraph.

For each $i \in \mathbb{Z}_m$, define $G_i := \{g_i \mid g \in G\}$. Analogous to classical Cayley digraphs, an $m$-Cayley digraph can be viewed as the digraph
\[
\Gamma = \mathrm{Cay}(G, T_{i,j} : i, j \in \mathbb{Z}_m)
\]
with:
\begin{itemize}
\item Vertex set: $G \times \mathbb{Z}_m = \bigcup_{i \in \mathbb{Z}_m} G_i$
\item Arc set: $\bigcup_{i,j \in \mathbb{Z}_m} \{(g_i, (tg)_j) \mid t \in T_{i,j}\}$
\end{itemize}
where $T_{i,j} \subseteq G$ for all $i,j \in \mathbb{Z}_m$.

When $T_{i,i} = \emptyset$ for all $i \in \mathbb{Z}_m$, we call $\Gamma$ an \emph{$m$-partite Cayley digraph}.
For any $g \in G$, the right multiplication map $R(g)$, defined by $R(g): x_i \mapsto (xg)_i$ for all $x_i \in G_i$ and $i \in \mathbb{Z}_m$, is an automorphism of $\Gamma$.$R(G)=\{R(g) \mid g \in G\}$ is isomorphic to $G$ and it is a semiregular group of automorphisms of $\Gamma$ with  $G_i$ as  orbits.

For a digraph $\Gamma$ and any vertex $x \in V(\Gamma)$, we recursively define the \emph{$k$-step out-neighborhood} of $x$ as follows:
\begin{itemize}
   \item $ \Gamma^{+0}(x) = \{x\}$,
  \item $\Gamma^{+1}(x) = \Gamma^+(x)$ is the set of out-neighbors of $x$ in $\Gamma$,
    \item $\Gamma^{+2}(x)$ represents the union of out-neighbors of all vertices in $\Gamma^+(x)$,
    \item  $\Gamma^{+k}(x) = \bigcup_{y \in \Gamma^{+(k-1)}(x)} \Gamma^+(y) \quad \text{for any integer } k \geq 1$.
 \end{itemize}

To end this section, we recall several important earlier results. First, Proposition 2.1 provides a powerful tool for our classification of finite groups admitting $m$-PDRs. Furthermore, by applying Propositions 2.2 and 2.3, we obtain a direct proof of Theorem 1.2.

\begin{Proposition}\cite[3.1]{du3}
Let $m$ be a positive integer at least $2$ and let $G$ be a finite group.
For any $i, j\in \mathbb{Z}_m$, let $T_{i,j}\subseteq G$ and let $\Gamma=Cay(G,T_{i,j}: i,j\in \mathbb{Z}_m)$ be a connected $m$-cayley digraph over $G$.
For $A=Aut(\Gamma)$, if $A$ fixes $G_i$ setwise for all $i\in \mathbb{Z}_m$ and there exist $u_0\in G_0,u_1\in G_1,\ldots,u_{m-1}\in G_{m-1}$ such that $A_{u_i}$ fixes $\Gamma^{+}(u_i)$ pointwise for all $i\in \mathbb{Z}_m$, then $A=R(G)$.
\end{Proposition}

\begin{Proposition}\cite[1.1]{gb}
Every finite simple group of order at least 5 has a ORR of out-valency 2.
\end{Proposition}

\begin{Proposition}\cite[2.4]{du4}
Let $G$ be a finite group and let $\mathrm{Cay}(G, R)$ be a DRR of $G$, where $1 \notin R \subset G$ and $|R| < |G|/2$. Then there exists a subset $L \subseteq G \setminus (R^{-1} \cup \{1\})$ with $|L| = |R|$ such that $\mathrm{Cay}(G, R \cup \{1\}, L \cup \{1\})$ forms a 2-PDR of $G$.
\end{Proposition}

\section{proof of Theorem}

\begin{Lemma}
Let \( G \) be an abelian group and $\Gamma=\mathrm{Cay}(G_i,T_{i,j}: i,j\in \mathbb{Z}_2)$ be a connected 2-partite Cayley digraph. If \( T_{0,1} = y T_{1,0} \), where \( y \) is an arbitrary element of \( G \), then the bijection $\tau$ defined as follows is an automorphism of $\Gamma$:
\[
\tau(g_0) = (yg)_1,\quad \tau(g_1) = g_0 \quad \text{for any element } g \in G.
\]
\end{Lemma}
\begin{proof}
\( \tau \in \text{Aut}(\Gamma) \) if and only if:
\[
(g_0, (sg)_1)^\tau = \left((g_0)^\tau, ((sg)_1)^\tau\right) = ((yg)_1, (sg)_0) \in \text{Arc}(\Gamma),
\]
and
\[
(g_1, (rg)_0)^\tau = \left((g_1)^\tau, ((rg)_0)^\tau\right) = (g_0, (yrg)_1) \in \text{Arc}(\Gamma),
\]
for all \( s \in T_{0,1} \) and \( r \in T_{1,0} \). This holds if and only if there exist \( r_1 \in T_{1,0} \) and \( s_1 \in T_{0,1} \) such that:
\[
sg = r_1 yg \quad \text{and} \quad yrg = s_1 g.
\]
This is evident because \( T_{0,1} = yT_{1,0} \).
\end{proof}

Lemma 3.1 establishes that for any finite abelian group $G$ and a 2-partite Cayley digraph $\Gamma = \mathrm{Cay}(G_i, T_{i,j}: i,j \in \mathbb{Z}_2)$, when $T_{0,1} = yT_{1,0}$ for some $y \in G$, we necessarily have $A > R(G)$, where $A$ is the full automorphism group of $\Gamma$.

The subsequent Lemmas 3.2 and 3.2 address the classification problem of finite cyclic groups $G = \langle x \rangle$ admitting $m$-PDRs of valency 3, with Lemma 3.2 handling the $m=2$ case and Lemma 3.3 covering all cases where $m \geq 3$.

\begin{Lemma}
Let \( G = \langle x \rangle \) be a finite cyclic group. Then, \( G \) admits 2-PDR of valency 3 if and only if the order of \( x \), denoted by \( o(x) \), satisfies \( o(x) \geq 5 \).
\end{Lemma}
\begin{proof}
Let $\Gamma=Cay(G_i,T_{i,j}: i,j\in \mathbb{Z}_2)$ be a connected 2-partite Cayley digraph and $A=Aut(\Gamma)$.
Since \( m = 2 \) and the valency is 3, we have \( o(x) \geq 3 \). If \( o(x)=3 \), then \( T_{0,1}=T_{1,0} \).
Therefore, it is evident that \( A > R(G) \) by Lemma 3.1.

\textbf{Case 1: $o(x)=4$}

Since $o(x)=4$, we have $G=\{1,x,x^2,x^3\}$.
We prove that there exists $0\leq i\leq 3$ such that $T_{1,0}=x^iT_{0,1}$.
Let's assume that $T_{1,0}\neq x^{i}T_{0,1}$ for all $0\leq i\leq 3$, then $T_{1,0}$ and $x^{i}T_{0,1}$ differ by exactly one element; otherwise $|T_{1,0}\cup x^{i}T_{0,1}|\geq 5>|G|$, which leads to a contradiction.

Suppose $T_{1,0}$ differs from $x^{i_1}T_{0,1}$ and $x^{i_2}T_{0,1}$ at elements $a$ and $b$ respectively where $0\leq i_1\neq i_2\leq 3$.
Then we must have $a=b$, since otherwise $|T_{1,0}\cup x^{i_1}T_{0,1}\cup x^{i_2}T_{0,1}|\geq 5>|G|$.
Since $G = \{1, x, x^2, x^3\}$, for a given $T_{1,0}$, we know from previous discussion that $x^i T_{0,1}$ has only three distinct possible forms for $0 \leq i \leq 3$.

To derive a contradiction, it suffices to show that $x^{i_1}T_{0,1}\neq x^{i_2}T_{0,1}$ for any $0\leq i_1<i_2\leq 3$. Without loss of generality, we prove $xT_{0,1}\neq x^2T_{0,1}$ (other cases follow similarly).

Assume $xT_{0,1}=x^2T_{0,1}$ and let $T_{0,1}=\{x^{j_1},x^{j_2},x^{j_3}\}$ with $0\leq j_1<j_2<j_3\leq 3$. Then:
\begin{align*}
xT_{0,1} &= \{x^{j_1+1},x^{j_2+1},x^{j_3+1}\} \\
x^2T_{0,1} &= \{x^{j_1+2},x^{j_2+2},x^{j_3+2}\}
\end{align*}

From $xT_{0,1}=x^2T_{0,1}$, we have $x^{j_1+1}\in x^2T_{0,1}$. Note that $0\leq j_1<j_2<j_3\leq 3$ implies $j_2-j_1,j_3-j_2\leq 2$ and $j_3-j_1\leq 3$.

Case analysis:
\begin{itemize}
\item If $x^{j_1+1}=x^{j_1+2}$, then $x=1$, contradicting $o(x)=4$.
\item If $x^{j_1+1}=x^{j_2+2}$, then $x^{j_2-j_1+1}=1$, which contradicts $o(x)=4$ and $j_2-j_1+1\leq 3$.
\item Therefore we must have $x^{j_1+1}=x^{j_3+2}$, i.e., $x^{j_3-j_1+1}=1$. Since $0\leq j_1<j_3\leq3$ and $j_3-j_1\leq 3$, this forces $j_3=3$ and $j_1=0$.
\end{itemize}

Similar arguments show $x^{j_2+1}=x^{j_1+2}=x^2$, yielding $j_2=1$. Thus $T_{0,1}=\{1,x,x^3\}$, but then clearly $xT_{0,1}\neq x^2T_{0,1}$, completing the proof.

Thus, when $o(x)=4$, there exists \( i \in \mathbb{Z}_4 \) such that \( T_{1,0} = x^i T_{0,1} \).
 By Lemma 3.1, we conclude that $A > R(G)$.

\textbf{Case 2: $o(x)= 5$ or $6$}

Let $T_{0,1}=\{1,x,x^2\}$ and $T_{1,0}=\{1,x,x^3\}$.

By using Mathematica or MAGMA \cite{MAG}, we obtain $A=R(G) \cong G$.
For the reader's convenience, we provide a rigorous mathematical proof for the case when $o(x)=5$, while the case for $o(x)=6$ can be proved similarly.

When $o(x)=5$, the graph $\Gamma$ is as shown in Fig 1.

\begin{figure}[H]
  \centering
  \includegraphics[width=1.0\linewidth]{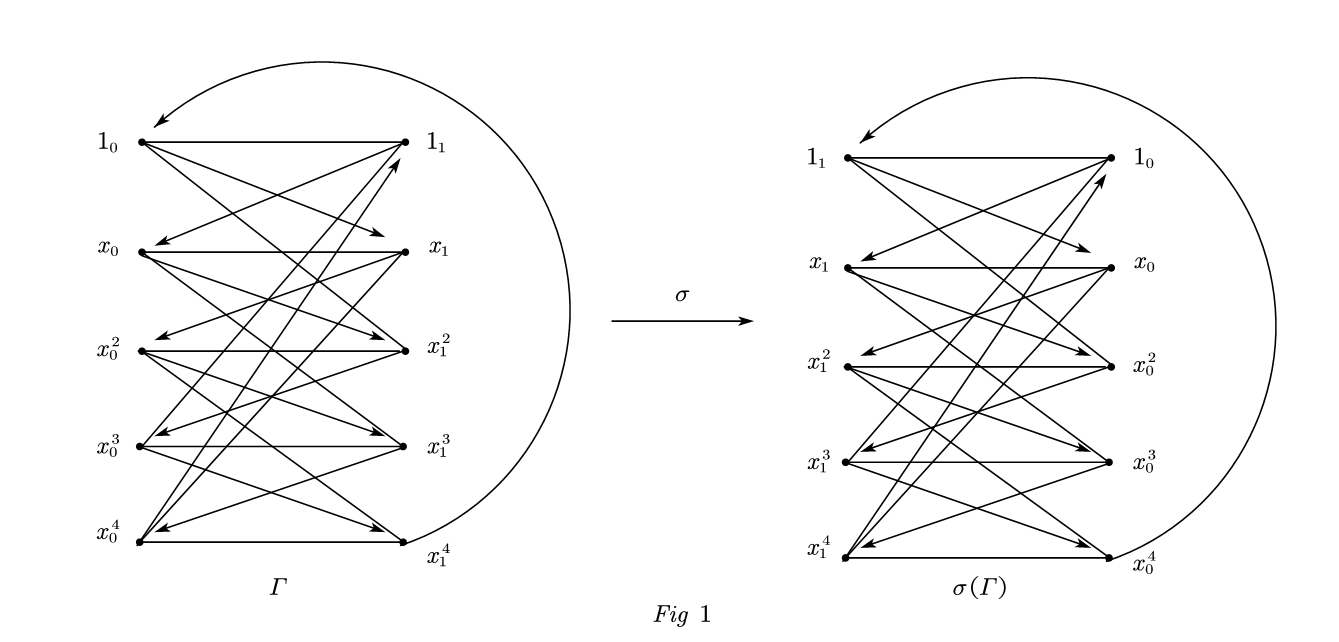}

\end{figure}

We first prove that $A$ fixes $G_0$ and $G_1$ setwise.
Suppose there exists $\sigma \in A$ and $g \in G$ such that $\sigma(g_0) \in G_1$. Without loss of generality, we assume $\sigma(1_0) = 1_1$ (this follows because if $\sigma(g_0) = t_1$ for some $g,t \in G$, then $(1_0)^{\rho_g\sigma\rho_{t^{-1}}} = 1_1$).

Note that $\Gamma$ contains a unique oriented cycle: $(1_0, x_1, x^2_0, x^3_1, x^4_0, 1_1, x_0, x^2_1, x^3_0, x^4_1)$.
It is readily observable that this cycle contains no undirected edges and passes through every vertex of $\Gamma$.
Therefore, when $\sigma(1_0) = 1_1$, we must have:
\begin{align*}
\sigma(x_1) = x_0,
\sigma(x^2_0) = x^2_1,
\sigma(x^3_1) = x^3_0,
\sigma(x^4_0) = x^4_1,
\sigma(1_1) = 1_0,
\sigma(x_0) = x_1,
\sigma(x^2_1) = x^2_0,
\sigma(x^3_0) = x^3_1,
\sigma(x^4_1) = x^4_0.
\end{align*}

The resulting $\sigma(\Gamma)$ is shown in Fig 1.
However, we observe that in $\Gamma$, $1_1$ is adjacent to $x^3_0$, while in $\sigma(\Gamma)$, $1_1$ is not adjacent to $x^3_0$. This contradicts $\sigma \in A$.
Therefore, $A$ must fix $G_0$ and $G_1$ setwise.
By Frattini argument, we have $A = R(G)A_{1_0}$.
Moreover, since $(1_0, x_1, x^2_0, x^3_1, x^4_0, 1_1, x_0, x^2_1, x^3_0, x^4_1)$ is the unique oriented cycle in $\Gamma$ and passes through every vertex of $\Gamma$, it follows that $A_{1_0} = 1$.
Consequently, we conclude that $A = R(G)$.

\begin{figure}[H]
  \centering
  \includegraphics[width=1.0\linewidth]{"fig2"}

\end{figure}

\textbf{Case 3: $o(x)\geq7$}

Let $T_{0,1}=\{1,x,x^2\}$ and $T_{1,0}=\{1,x,x^3\}$, then we have:

\[
\begin{aligned}
&\Gamma^{+}(1_0) = \{1_1,x_1,x^2_1\}, \quad \Gamma^{+}(1_1) = \{1_0,x_0,x^3_0\} \\
&\Gamma^{+2}(1_0) = \{1_0,x_0,x^2_0,x^3_0,x^4_0,x^5_0\}, \quad \Gamma^{+2}(1_1) = \{1_1,x_1,x^2_1,x^3_1,x^4_1,x^5_1\} \\
&\Gamma^{+3}(1_0) = \{1_1,x_1,x^2_1,x^3_1,x^4_1,x^5_1,x^6_1,x^7_1\}, \quad \Gamma^{+3}(1_1) = \{1_0,x_0,x^2_0,x^3_0,x^4_0,x^5_0,x^6_0,x^7_0,x^8_0\} \\
\end{aligned}
\]

First, we demonstrate that $A$ fixes $G_0$ and $G_1$ setwise.

When $o(x) \geq 9$, we clearly have $|\Gamma^{+3}(1_1)| = 9 > |\Gamma^{+3}(1_0)| = 8$, which implies that $A$ fixes $G_0$ and $G_1$ setwise.

For the cases where $o(x)=7$ or $8$, as shown in Fig 2, $H_i$ represents the sub-digraphs induced by $\bigcup_{0 \leq k \leq 2} \Gamma^{+k}(1_i)$ (for clarity of observation, we have omitted the arcs from $\Gamma^{+1}(1_i)$ to $\Gamma^{+2}(1_i)$), where $i = 0$ or $1$.We observe that in $\Gamma^{+2}(1_0)$, there exists a vertex $x_0$ with out-degree $2$ in $\Gamma^{+}(1_0)$, whereas no such vertex exists in $\Gamma^{+2}(1_1)$. Therefore, $A$ fixes $G_0$ and $G_1$ setwise.
In summary, we have established that $A$ fixes both $G_0$ and $G_1$ setwise.

We now prove that $A_{1_0} = 1$.
By the Frattini argument, we have:
\[ A = R(G)A_{1_0} = R(G)A_{x_0} = R(G)A_{1_1} = R(G)A_{x_1}, \]
which implies:
\[ |A_{1_0}| = |A_{x_0}| = |A_{1_1}| = |A_{x_1}| = |A/R(G)|. \]

We observe that in the subgraph $[\Gamma^{+}(1_0)]$, the vertices $1_0$ and $1_1$ are uniquely connected by an undirected edge.
Thus, $A_{1_0}$ fixes $\{1_1\}$ and $\{x_1,x^2_1\}$ setwise.
Moreover, since $x_1$ and $x^2_1$ have in-valency $1$ and $2$ in $\Gamma^{+2}(1_0)$, respectively, $A_{1_0}$ fixes $x_1$ and $x^2_1$. The vertex $x_0$ is the unique in-neighbor of $x_1$ in $\Gamma^{+2}(1_0)$, so $A_{1_0}$ also fixes $x_0$. Consequently:
\[ A_{1_0} = A_{x_0} = A_{1_1} = A_{x_1}. \]

We now prove that $A_{x_0}=A_{x_0^2}$ using  $A_{1_0}=A_{x_0}$.
For $\sigma \in A_{x_0}$, since
\[ 1_0^{\rho_x\sigma\rho_{x^{-1}}} = x_0^{\sigma\rho_{x^{-1}}} = x_0^{\rho_{x^{-1}}} = 1_0, \]
we have $\rho_x\sigma\rho_{x^{-1}} \in A_{1_0}=A_{x_0}$.
This implies that $x_0^{\rho_x\sigma\rho_{x^{-1}}} = x_0$, which is equivalent to $(x_0^2)^\sigma = x_0^2$.
Therefore, $\sigma \in A_{x_0^2}$, and consequently $A_{x_0}=A_{x_0^2}$ by  $|A_{x_0}|=|A_{x_0^2}|$.

From the preceding discussion and by considering the conjugate action of $R(G)$ on $G_i$, we have:
\[ A_{x_0} = A_{x^2_0}= A_{x^3_0}=\dots= A_{1_0} =A_{x_1}= A_{x^2_1} =A_{x^3_1}= \dots = A_{1_1}. \]

Since $G=\langle x \rangle$ is a cyclic group, it follows that $A_{1_0}=1$, and thus $A=R(G) \cong G$.

\end{proof}

\begin{Lemma}
Let $G = \langle x \rangle$ with $o(x) \geq 2$ and $m \geq 3$. Then $G$ admits a $m$-PDR of valency 3, except when \( o(x) = 2 \) and \( m=3 \).
\end{Lemma}

\begin{proof}
Let $\Gamma=Cay(G_i,T_{i,j}: i,j\in \mathbb{Z}_m)$ be a connected $m$-partite Cayley digraph and $A=Aut(\Gamma)$.

\textbf{Case 1:  $o(x)=2$ and $m=3$ }

When \( m = 3 \), the valency condition  ($\Gamma$ being 3-regular) implies that either \( T_{0,1} = G \) or \( T_{0,2} = G \) must hold.
Without loss of generality, assume $T_{0,2}=G$. We shall prove that $T_{0,2}=T_{2,1}=T_{1,0}=G$.

If $T_{2,0}=G$, then $T_{1,0}\neq G$ (otherwise vertex $G_0$ would have in-valency 4, leading to a contradiction). Consequently, $T_{1,2}=G$, but this would result in $G_2$ having in-valency 4, which is again a contradiction. Therefore, $T_{2,0}\neq G$, and by the valency 3 condition, we must have $T_{2,1}=G$. Similarly, we obtain $T_{1,0}=G$.

Therefore, the valency condition implies $\lvert T_{0,1} \rvert = \lvert T_{1,2} \rvert = \lvert T_{2,0} \rvert = 1$. 
However, for any assignment of $T_{0,1},T_{1,2},T_{2,0}$ to either $1$ or $x$, our Mathematica computations show that $|A|=6>|G|$.

\textbf{Case 2:  $o(x)=2$ and $m=4$ }

\begin{figure}[H]
  \centering
  \includegraphics[width=0.7\linewidth]{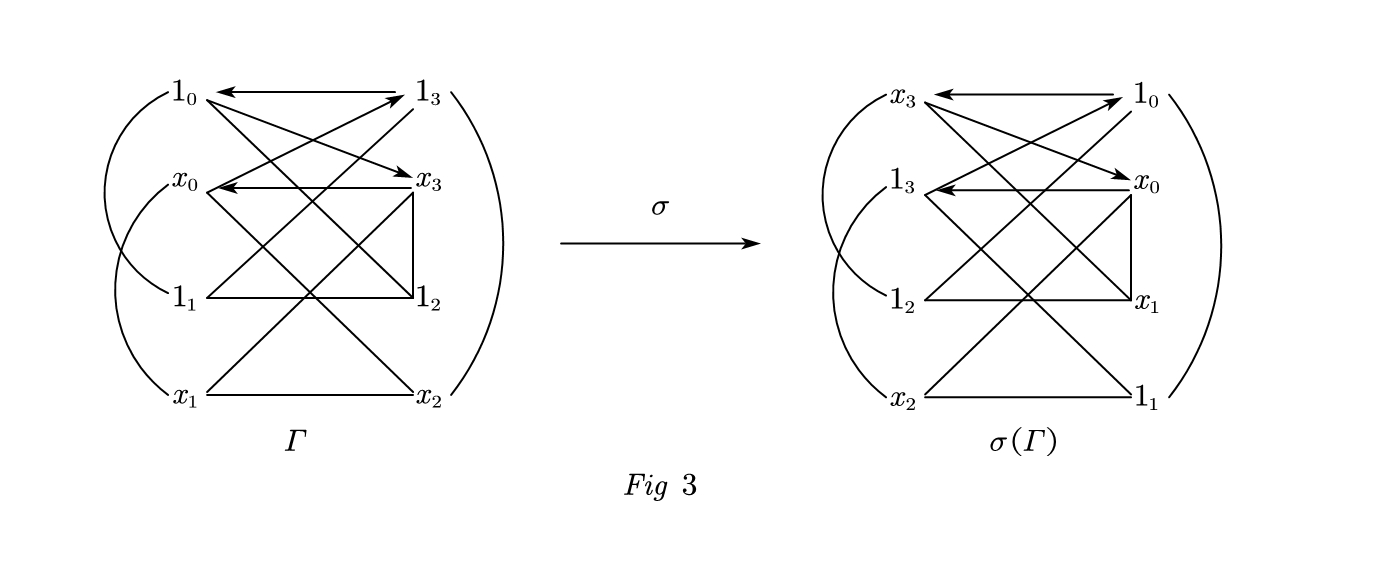}

\end{figure}

We set \(T_{0,1} = T_{0,2}=\{1\}\), \(T_{0,3}=x\); \(T_{1,0}=T_{1,2}=T_{1,3}=\{1\}\); \(T_{2,0}=T_{2,1}=\{1\}\), \(T_{2,3}=\{x\}\); \(T_{3,0}=T_{3,1}=\{1\}\), \(T_{3,2}=\{x\}\). Then \(\Gamma\) is shown in Fig 3.
Next, we will prove that \(A = R(G)\).
Notice that, except for the lack of an undirected edge between \(G_0\) and \(G_3\), there is an undirected edge between any \(G_i\) and \(G_j\) for $i\neq j \in \mathbb{Z}_4$.
So \(A\) fixes \(G_0\cup G_3\) setwise.
If there exists \(\sigma\in A\) such that \(\sigma(g_0)\in G_3\), without loss of generality, let \(\sigma(1_0)=1_3\).
Notice that, in \([G_0\cup G_3]\), \((1_0,x_3,x_0,1_3)\) is the only oriented cycle.
So we must have \(\sigma(x_3)=1_0\), \(\sigma(x_0)=x_3\), \(\sigma(1_3)=x_0\).
Moreover, in \(G_1\cup G_2\), the only vertex adjacent to both \(1_0\) and \(1_3\) is \(1_1\), the only vertex adjacent to both \(x_0\) and \(x_3\) is \(x_1\), the only vertex adjacent to both \(1_0\) and \(x_3\) is \(1_2\), and the only vertex adjacent to both \(x_0\) and \(1_3\) is \(x_2\) (this means that if we fix \(G_0\cup G_3\) pointwise, then we fix V(\(\Gamma\)) pointwise).
So we have \(\sigma(1_1)=x_2\), \(\sigma(x_1)=1_2\), \(\sigma(1_2)=1_1\), \(\sigma(x_2)=x_1\), that is, \(\sigma(\Gamma)\) is shown in Fig 3.

However, we notice that in \(\Gamma\), \(1_1\) and \(x_2\) are non-adjacent, but in \(\sigma(\Gamma)\), \(1_1\) and \(x_2\) are adjacent, which contradicts \(\sigma\in A\).
So \(A\) fixes \(G_0\) and \(G_3\) setwise. By the Frattini argument, we know that \(A = R(G)A_{1_0}\).
Next, we prove that \(A_{1_0}=1\).
From the previous discussion, we know that we only need to prove that \(A_{1_0}\) fixes \(\{1_0,x_0,1_3,x_3\}\) pointwise.
This is obvious because \((1_0,x_3,x_0,1_3)\) is the only oriented cycle in \([G_0\cup G_3]\).
So \(A_{1_0}=1\), that is, \(A = R(G)\).

\textbf{Case 3: $o(x)=2$ and $m\geq5$}

Define the subsets of $G$ as follows:
\[
T_{i,i-1}=T_{i,i+1}=\{1\},\quad i\in \mathbb{Z}_m
\]
\[
T_{i,i+2}=\{1\},\quad T_{0,2}=\{x\}
,\quad i\in \mathbb{Z}_m\setminus\{0\}
\]
with all other $T_{i,j}=\emptyset$, where $i,j\in \mathbb{Z}_m$.

The graph \(\Gamma\) is shown in Fig 4.
\begin{figure}[H]
  \centering
  \includegraphics[width=0.9\linewidth]{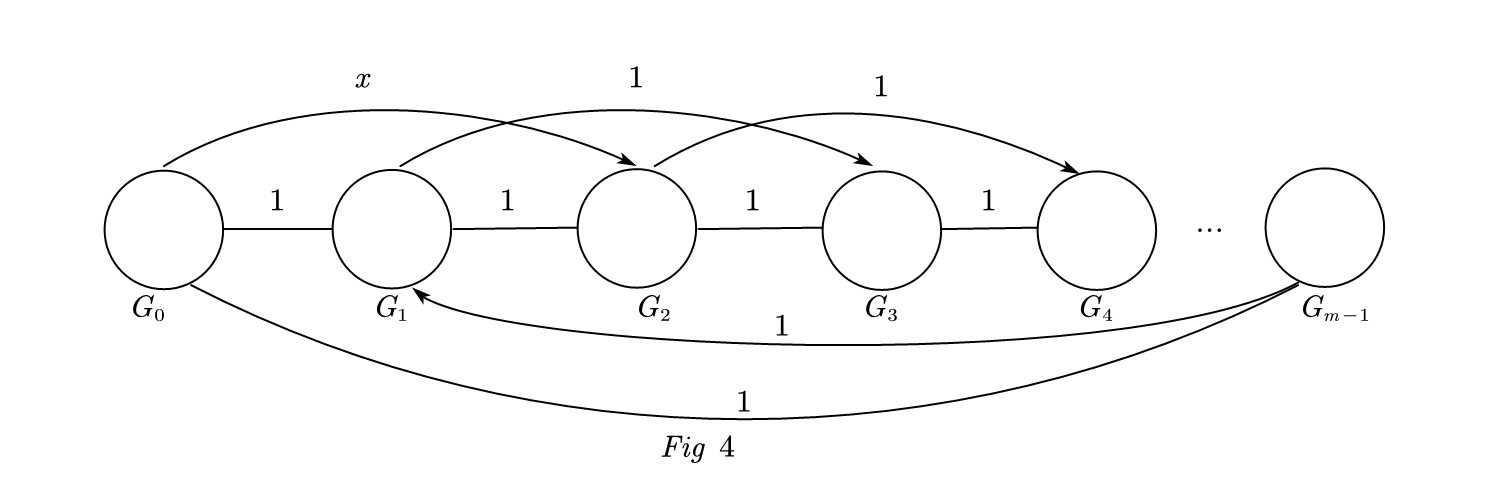}

\end{figure}
Observe that $(g_i,g_{i+2},g_{i+1})$ forms a directed 3-cycle for all $g \in G$, with the exception of the special case $(g_0,g_2,g_1)$.
This implies that $A$ fixes $G_0\cup G_1\cup G_2$ setwise.
Considering the subgraph $[G_0\cup G_1\cup G_2]$, we note that vertices in $G_1$ have two undirected edges: $g_0\sim g_1$ and $g_1\sim g_2$, but each vertex in $G_0$ and $G_2$ is incident to exactly one undirected edge.
Consequently, $A$ fixes $G_1$ setwise.

Furthermore, in $[G_0\cup G_1\cup G_2]$, vertices in $G_0$ have out-valency 2 while those in $G_2$ have out-valency 1, which implies that $A$ fixes $G_0$, $G_1$, and $G_2$ setwise. Since $T_{1,3}=\{1\}$ and $T_{1,i}=\emptyset$ for $i\neq 0,2,3$, $A$ fixes $G_3$ setwise.
Similarly, since \( T_{2,i} = \emptyset \) for \( i \neq 1,3,4 \), we conclude that \( A \) fixes \( G_4 \) setwise.
Using the same method, we can show that \( A \) fixes all \( G_i \) setwise for \( 0 \leq i \leq m-1 \).

Since $\Gamma^{+}(1_i)$ lies in distinct parts and $A$ fixes each $G_i$ setwise for $0\leq i\leq m-1$, we have $A_{1_i}$ fixes each vertex in $\Gamma^{+}(1_i)$  and Proposition 2.1 implies that $A=R(G)\cong G$.

\textbf{Case 4: $o(x)\geq3$ and $m\geq3$}

Define the subsets of $G$:
\[
T_{i,i+1}=\{1,x\},\quad T_{j,j-1}=\{1\}\quad (i,j\in \mathbb{Z}_m,\; j\neq1),\quad T_{1,0}=\{x\}
\]
with all other $T_{i,j}=\emptyset$.

The graph \(\Gamma\) is shown in Fig 5.
\begin{figure}[H]
  \centering
  \includegraphics[width=0.9\linewidth]{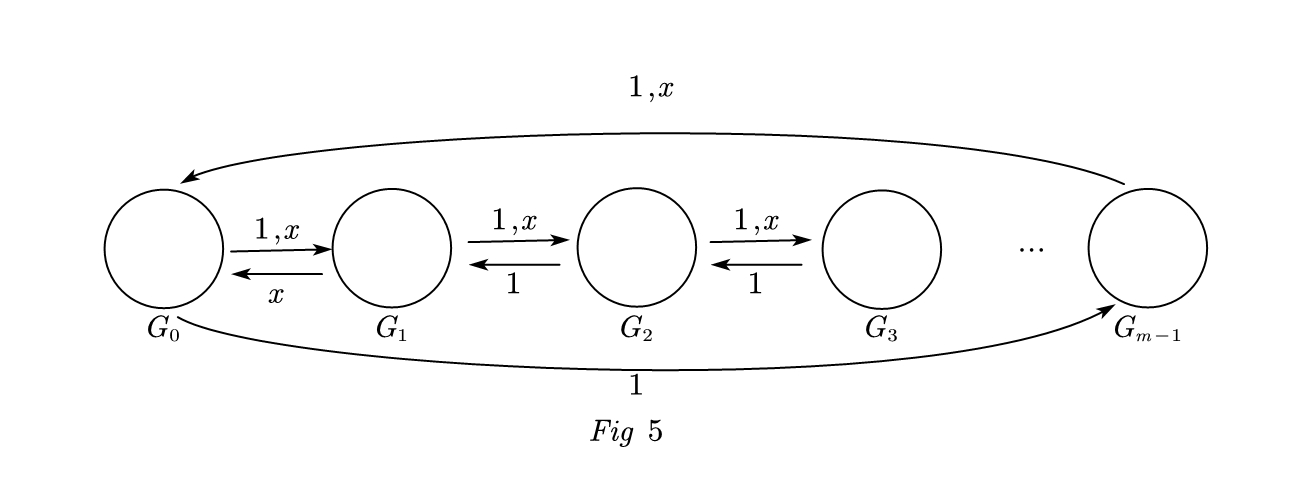}

\end{figure}
The edge $g_i\sim g_{i+1}$ is undirected for $i\neq0$, but the subgraph $[G_0\cup G_1]$ contains no undirected edges.
Therefore, $A$ fixes $G_0\cup G_1$ setwise.
In subgraph $[G_0\cup G_1]$, vertices in $G_0$ have out-valency 2 while those in $G_1$ have out-valency 1, so $A$ fixes both $G_0$ and $G_1$ setwise.

Since $T_{1,i}=0$ for $i\neq0,2$, $A$ fixes $G_2$.
Similarly, we can obtain that $A$ fixes all $G_i$ setwise for $0\leq i\leq m-1$.
Then we have $A=R(G)A_{g_i}$ by Frattini argument where $g_i\in G_i$.

From the fact that $A_{1_0}$ fixes $1_{m-1}$ by  $T_{0,m-1}=\{1\}$, we obtain $A_{1_0}=A_{1_{m-1}}$ by Frattini argument.
Since $A_{1_{m-1}}$ fixes $1_{m-2}$ by $T_{m-1,m-2}=\{1\}$, it follows that $A_{1_0}$ fixes $1_{m-2}$, and thus $A_{1_0}=A_{1_{m-1}}=A_{1_{m-2}}$. Continuing this process, we conclude that $A_{1_0}=A_{1_1}=A_{1_2}=\cdots=A_{1_{m-1}}$.

Moreover, since $A_{1_1}$ fixes $x_0$ by $T_{1,0}=\{x\}$, we have $A_{1_1}=A_{x_0}$.
Repeating this argument yields $A_{1_1}=A_{x_0}=A_{x_1}=A_{x_2}=\cdots=A_{x_{m-1}}$.

By the conjugate action of $R(G)$ on $G_i$, we establish that:
\[ A_{x_0} = A_{x^2_0}= \dots= A_{1_0} =A_{x_1}= A_{x^2_1} = \dots = A_{1_1}=\dots=A_{x_{m-1}} = A_{x^2_{m-1}}= \dots= A_{1_{m-1}}. \]

Since $G$ is a cyclic group, $A_{1_0}=1$, and by the Frattini argument, we conclude that $A=R(G)\cong G$.
\end{proof}

The conclusion of Theorem 1.1 follows directly from Lemmas 2.3 and 2.4.
We now consider the case where $G = \langle x, y \rangle$ is a two-generated group.

\begin{proof}[\textbf{Proof of Theorems  1.2}]
Since $G=\langle x,y\rangle$ is a finite simple group, we have $|G|\geq 6$. By Proposition 2.2, there exists a subset $R\subseteq G$ such that $\mathrm{Cay}(G,R)$ is an ORR, where $|R|=2<|G|/2$. Because $\mathrm{Cay}(G,R)$ is connected and $|R|=2$ (an ORR is clearly a DRR), we have $1\notin R$ ($\mathrm{Cay}(G, R)$ is connected if and only if $G=\langle R\rangle$). Therefore, by Proposition 2.3, there exists a subset $L\subseteq G$ such that $\mathrm{Cay}(G,R\cup\{1\},L\cup\{1\})$ is a 2-PDR of valency 3, where $L\subseteq G\setminus(R^{-1}\cup\{1\})$ and $|L|=|R|=2$.

\end{proof}

\begin{proof}[\textbf{Proof of Theorems  1.3}]

\begin{figure}[H]
  \centering
  \includegraphics[width=0.9\linewidth]{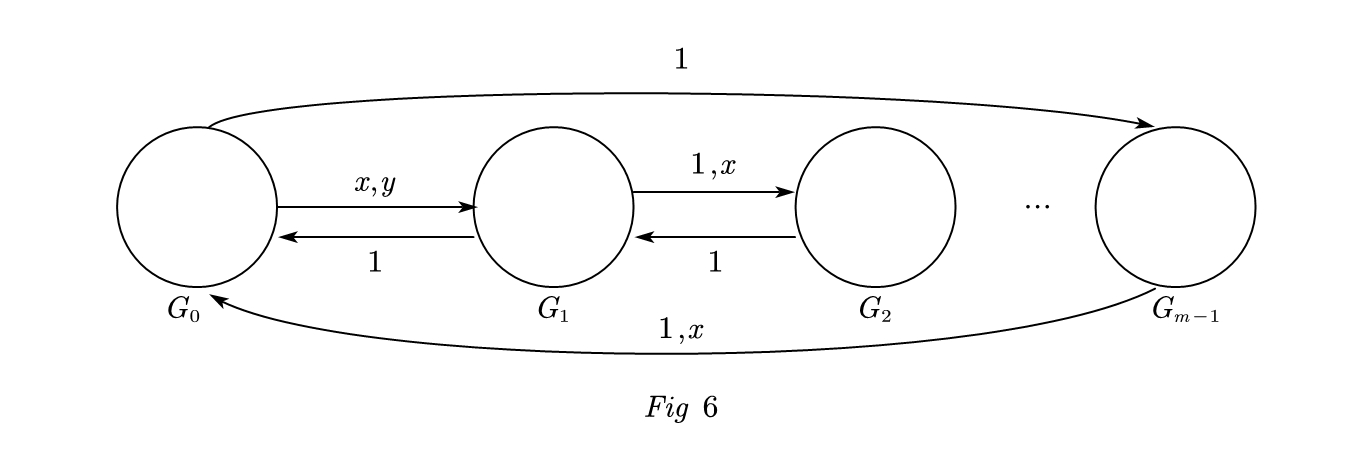}

\end{figure}
Define the connection sets as follows:
\begin{align*}
T_{0,1} &= \{x,y\}, \\
T_{i,i+1} &= \{1,x\} \text{ for } i \in \mathbb{Z}_m \text{ with } i \neq 0, \\
T_{j,j-1} &= \{1\} \text{ for } j \in \mathbb{Z}_m, \\
T_{i,j} &= \emptyset \text{ for all other cases.}
\end{align*}

Thus we obtain:
\begin{align*}
T_{0,1} \cap T^{-1}_{1,0} &= \emptyset, \\
T_{i,i+1} \cap T^{-1}_{i+1,i} &= \{1\} \text{ for all } i\in \mathbb{Z}_m\setminus \{0\}.
\end{align*}

Let $\Gamma := \mathrm{Cay}(G, T_{i,j} : i,j \in \mathbb{Z}_m)$, then $\Gamma$ is as shown in Fig 6, and let $A = \mathrm{Aut}(\Gamma)$.
 Clearly, $\Gamma$ is a 3-regular $m$-partite Cayley digraph of $G$. To complete our proof, it suffices to establish that $A \cong G$.

Observe that:
\begin{itemize}
\item Vertices in $G_0$ and $G_1$ are each incident to exactly one undirected edge.
\item Vertices in $G_i$ (for $i \neq 0,1$) are each incident to two undirected edges.
\end{itemize}

This implies that $A$ fixes $G_0 \cup G_1$ setwise. Furthermore, in the subgraph $[G_0 \cup G_1]$:
\begin{itemize}
\item Vertices in $G_0$ have out-valency 2.
\item Vertices in $G_1$ have out-valency 1.
\end{itemize}
Therefore, $A$ fixes both $G_0$ and $G_1$ setwise.
Note that since $T_{1,i} = \emptyset$ for all $i \neq 0, 2$, the group $A$ fixes $G_2$ setwise.
By similar arguments, we deduce that $A$ fixes $G_i$ setwise from $T_{i,j}=\emptyset$ for all $0 \leq i \leq m-1$ and $j\neq i-1, i+1$.

We now prove that $A_{1_i}$ fixes $\Gamma^{+}(1_i)$ pointwise.
Since $T_{i,i-1} = 1$, we have the stabilizer chain:
\[ A_{1_0} = A_{1_1} = A_{1_2} = \cdots = A_{1_{m-1}}. \]

Observe that $\Gamma^{+}(1_{m-1}) = \{1_0, x_0, 1_{m-2}\}$, which implies that $A_{1_0}$ fixes $x_0$ and $A_{1_0}=A_{x_0}$.
From $T_{1,0} = \{1\}$, we conclude that $A_{x_0}$ fixes $x_1$ and $A_{x_0}=A_{x_1}$.
Consequently, $A_{1_0}$ fixes $\Gamma^{+}(1_0) = \{1_{m-1}, x_1, y_1\}$ pointwise.

Similarly, we can prove that $A_{1_i}$ fixes $\Gamma^{+}(1_i)$ pointwise for all $1 \leq i \leq m-1$. Applying Proposition 2.1, we finally obtain $A = R(G)$.

\end{proof}

\section{Acknowledgments}

We gratefully acknowledge the support of the Graduate Innovation Program of China University of Mining and Technology (2025WLKXJ146), the Fundamental Research Funds for the Central Universities, and the Postgraduate Research and Practice Innovation Program of Jiangsu Province (KYCX25\_2858) for this work.

\end{document}